\documentclass[11pt]{amsart}
\begin{document}

\title{Injective Envelopes of $C^*$-algebras as Operator Modules}

\author[M.~Frank]{Michael Frank*}
\address{Universit\"at Leipzig, Mathematisches Institut,
D-04109 Leipzig, Germany}
\email{frank@mathematik.uni-leipzig.de}
\thanks{*Research supported in part by a grant from the NSF}
\author[V.~I.~Paulsen]{Vern I.~Paulsen*}
\address{Dept.~of Mathematics, University of Houston, Texas  77204-3476,
U.S.A.}
\email{vern@math.uh.edu}

\begin{abstract}
In this paper we give some characterizations of M.~Hamana's injective
envelope $I(A)$ of a $C^*$-algebra $A$ in the setting of operator spaces and completely
bounded maps.  These characterizations lead to simplifications and
generalizations of some known results concerning completely bounded
projections onto $C^*$-algebras.
We prove that $I(A)$ is rigid for completely bounded $A$-module maps.
This rigidity yields a natural representation of many kinds of multipliers as 
multiplications by elements of $I(A)$.  In particular, we prove that the($n$ times 
iterated) local multiplier algebra of $A$ embeds into $I(A)$.
\end{abstract}
\maketitle


\baselineskip=20pt

\newtheorem{thm}{Theorem}[section]
\newtheorem{prop}[thm]{Proposition}
\newtheorem{lem}[thm]{Lemma}
\newtheorem{cor}[thm]{Corollary}
\theoremstyle{definition}
\newtheorem{dfn}{Definition}
\theoremstyle{remark}
\newtheorem*{rmk}{Remark}
\newtheorem*{prob}{Problem}
\newtheorem*{exam}{Example}
\maketitle

\newcommand{\Cal}{\mathcal}
\newcommand{\rE}{\tilde{E}}
\def\lam{\lambda}
\def\alg{\text{ alg}}
\def\ball{\text{ ball}}
\def\dist{\text{ dist}}
\def\diag{\text{ diag}}
\def\var{\varphi}
\def\span{\text{ span}}
\def\Im{\text{Im }}
\def\D{\mathbb D}
\def\id{\text{ id}}
\newcommand{\n}[1]{\left\vert#1\right\vert}   
\def\bv \bign#1{\bigl\vert#1\bigr\vert} 
\newcommand{\V}[1]{\left\Vert#1\right\Vert}   
\def \p#1{\left(#1\right)}         
\def\ra{\rightarrow}
\def \a#1{\langle#1\rangle}
\def \supnorm#1{\N{#1}_{\infty}}        
\def \onenorm#1{\N{#1}_{1}}             
\def \Xa{\bold X}
\def \R{\mathbb R}
\def\cO{\Cal O}
\def\cX{\Cal X}
\def\cM{\Cal M}
\def\cZ{\Cal Z}
\def\cA{\Cal A}
\def\cN{\Cal N}
\def\cB{\Cal B}
\def\cU{\Cal U}
\def\cH{\Cal H}
\def\cJ{\Cal J}
\def\cI{\Cal I}
\def\cR{\Cal R}
\def\cS{\Cal S}
\def\cK{\Cal K}
\def\be{\beta}
\def\ep{\epsilon}
\def\al{\alpha}
\def\raw{\overset\rightarrow \to w}

\def\ralam{\overset\rightarrow\to  \lambda}
\def\ramu{\overset\rightarrow\to \mu}
\def\raz{\overset\rightarrow \to z}
\def\rax{\overset\rightarrow \to x}
\def\ray{\overset\rightarrow \to y}
\def\raal{\overset\rightarrow\to \alpha}
\def\rabe{\overset\rightarrow\to\beta}
\def\ga{\gamma}
\def\raga{\overset\rightarrow\to\gamma}
\def \Rn{\R^n}                   
\def \Z{\mathbb Z}
\def\C{\mathbb C}              
\def\cS{\Cal S}
\def \cC{\Cal C}
\def\cD{\Cal D}
\def \F{\mathbb F}
\def \N{\mathbb N}
\def\cP{\Cal P}
\def\cF{\Cal F}
\def \Sa{\mathbb S}
\def \MIN{\text{MIN}}
\def \MAX{\text{MAX}}
\def\raw{\overset\ra\to w}
\def\FO{\overset\frown\to\otimes}
\def\SO{\underset\smile\to\otimes}
\def\pp{\prime\prime}
\def\Bar{\overline}


\section{Introduction} Let $A$ denote a unital $C^*$-algebra.
M.~Hamana \cite{Ham,Ham2,Ham3} introduced the injective envelope of
$A, \ I(A),$ as a ``minimal" injective operator system containing $A$
and established various characterizations and properties of $I(A)$ in
the setting of completely positive mappings and operator systems.

In recent years attention has shifted from completely positive maps
and operator systems to completely bounded maps, operator algebras and
ope\-rator spaces. In particular a theory has evolved of operator spaces
that
are completely contractive as modules over operator algebras.  See
for example, \cite{BMP,Pa}.

This theory gives a new categorical framework where one can examine
injective
envelopes.  While other author's have pursued this viewpoint they have
generally defined injectivity and rigidity in terms of completely
contractive
maps.  For example, defining injectivity by requiring completely contractive
maps to have completely contractive extensions. This is equivalent to
requiring that completely bounded maps have completely bounded extensions of
the same completely bounded norm. Since unital completely contractive maps
on $C^*$-algebras are completely positive this approach ge\-nerally reduces 
to M.~Hamana's results in the $C^*$-algebra setting.

Our approach is different in that we are interested in a setting
where our objects are $A$-modules and injectivity is defined by
requiring that completely bounded $A$-module maps have completely bounded
$A$-module extensions, but not necessarily of the same norm. We
show for example that, as well as being a minimal injective operator system
containing $A,$ that $I(A)$ is in a certain sense the ``minimal" injective
left operator $A$-module containing $A.$  This immutability of the
``injective hull" of $C^*$-algebras under change of ca\-tegory has some
immediate applications to completely bounded projections and multiplier
algebras.

In this paper our primary focus is on the new properties of $I(A)$
and their applications.  For this reason we take the quickest approach,
which is to restrict our attention to the unital case and use M.~Hamana's 
results to deduce these new properties by our off-diagonalization trick.

Many of our results do carry over to the case of a non-unital $C^*$-algebra
$B$ by the simple device of adjoining a unity to $B$, letting $A$ denote
this unital $C^*$-algebra and observing that any $B$-modules are
automatically $A$-modules. For a greater development of the non-unital case
we refer the reader to the subsequent paper \cite{BP}.


\section{Mapping Properties of $I(A)$}

Throughout this section $A$ will denote a unital $C^*$-algebra and $I(A)$
will denote its injective envelope as defined in \cite[Def.~2.1,
Th.~4.1]{Ham}.
We assume that the reader is familiar with the definitions and elementary
properties of completely bounded and completely positive maps as presented
in
\cite{Pa} or \cite{Pi}.

One of M.~Hamana's fundamental results about $I(A)$ was his rigidity
theorem.
This theorem says that if $\var :I(A) \ra I(A)$ is completely positive with
$\var (a) = a$ for all $a$ in $A$ then $\var (x) = x$ for all $x$ in $I(A).$

A direct analog of this result is false for general completely bounded maps.
If $A \neq I(A),$ then there exists a non-zero bounded linear functional
$f: I(A) \ra \C$ with $f(A) = \{ 0 \}.$ Defining the map $\var : I(A) \ra
I(A)$
via $\var (x) = x + f(x) 1$ yields a completely bounded map with $\var (a) =
a$
for all $a$ in $A,$ but $\var (x) \neq x$ for all $x $ in $I(A).$
However, if one recalls that completely positive maps that fix $A$
are automatically $A$-bimodule maps \cite{Choi}, \cite[Exercise 4.3]{Pa}
then one is led to the appropriate generalization of M.~Hamana's rigidity.
Surprisingly one does not need bimodules, only left or right $A$-modules as
the following results show.

\begin{thm}  \label{thm-1}
   Let $E \subseteq I(A)$ be a subspace such that $A E \subseteq E$
   (respectively, $E A \subseteq E$) and let $\var : E \ra I(A)$ be a
   completely bounded left(resp., right) $A$-module map.
   Then there exists an
   element $y$ in $I(A)$ such that $\var$ is right (resp., left)
   multiplication by $y,$ i.e., $\var (e) = e y$ ($\var (e) = ye$) for
   all $e$ in $E$ and $\V{y} = \V{\var}_{cb}.$  When $AEA \subseteq E$
   and $\var$ is a bimodule map, then $y$ may be taken in the center of
   $I(A).$

   In particular, $\var$ extends to a completely bounded, left(resp., right)
   $A$-module map $\psi: I(A) \to I(A)$ such that $\psi|_E = \var$ and
   $\| \var \|_{cb} = \| \psi \|_{cb}$. If $E \subseteq I(A)$ contains an
   invertible element of $I(A)$, then the element $y \in I(A)$ and
   consequently the extension $\psi$ are unique.
\end{thm}

\begin{proof}
It will suffice to assume that $\V{\var}_{cb} \leq 1.$  Let $\cS \subseteq
M_2 (I(A))$ be defined as
  $$
    \cS = \left\{ \begin{pmatrix}
    a & e \\
    f^* & \lam
    \end{pmatrix} : a \in A, e, f, \in E, \lam \in \C \right\}
  $$
and $\Phi : \cS \ra M_2 (I(A))$ by
  $$
    \Phi \left( \, \begin{pmatrix}
    a & e \\
    f^* & \lam \end{pmatrix} \, \right) = \begin{pmatrix}
    a & \var (e) \\
    \var (f)^* & \lam
    \end{pmatrix}.
  $$
Arguing as in \cite{Pa} or \cite{Suen} one sees that $\Phi$ is completely
positive and hence can be extended to a completely positive
map on all of $M_2 (I(A))$ which we still denote by $\Phi.$
Using the fact that $\Phi$ fixes $A \oplus \C = \left\{\begin{pmatrix}
a & 0 \\
0 & \lam
\end{pmatrix} : a \in A, \lam \in \C\right\}$ and again arguing as in
\cite{Suen}, we see
that there exists $\var_i : I(A) \ra I(A),         \linebreak[4]
i = 1, 2, 3, 4$ such that
  $$
    \Phi \left( \, \begin{pmatrix}
    x_1 & x_2 \\
    x_3 & x_4
    \end{pmatrix} \, \right) =
    \begin{pmatrix}
    \var_1 (x_1) & \var_2 (x_2)\\
    \var_3 (x_3) & \var_4 (x_4)
    \end{pmatrix}.
  $$
Clearly, $\var_1$ and $\var_4$ must be completely positive and $\var_2$
extends $\var.$

Since $\var_1 (a) = a$ for all $a$ in $A,$ by M.~Hamana's rigidity result
$\var_1 (x) = x$ for all $x$ in $I(A).$  Thus, $\Phi$ fixes the
$C^*$-subalgebra, $I(A) \oplus \C$ and so by     \linebreak[4]
\cite{Choi} (see also \cite{Pa}) $\Phi$ must be a bimodule map over this
algebra. Thus,   \linebreak[4]
$\begin{pmatrix} 0 & \var_2 (x)\\
0 & 0 \end{pmatrix} = \Phi \left( \! \begin{pmatrix} 0 & x\\
0 & 0 \end{pmatrix} \! \right) = \Phi \left( \! \begin{pmatrix}
x & 0 \\
0 & 0 \end{pmatrix} \! \begin{pmatrix}
0 & 1 \\
0 & 0 \end{pmatrix} \! \right)=
\begin{pmatrix}
x & 0\\
0 & 0
\end{pmatrix} \! \begin{pmatrix}
0 & \var_2 (1)\\
0 & 0
\end{pmatrix}$
and we have that $\var_2 (x) = x \cdot \var_2 (1).$
Finally, $\V{\var}_{cb} = \V{\var_2}_{cb} = \V{\var_2 (1)}.$

The proof for right $A$-module maps is similar.  For the bimodule case
let $\cS = \left\{\begin{pmatrix}
a & e\\
f^* & b
\end{pmatrix} : a, b \in A, e, f \in E \right\}$
and deduce that $\var_2$ is an $I(A)$-bimodule map.
If $E$ contains an invertible element e, then $y = e^{-1} \var (e)$
and so $y$ is unique.
\end{proof}

In particular, the above results show that
every completely bounded left (resp., right or bi-)
$A$-module map of $A$ into $I(A)$ admits a unique extension to a
completely bounded left (resp., right or bi-) $A$-module map of $I(A)$
into itself and this extension has the same completely bounded norm.

\begin{cor}[Rigidity]  \label{cor-rigid}
   Let $A$ be a unital $C^*$-algebra and $I(A)$ be its injective envelope
   $C^*$-algebra. Let $E$ be a subspace of $I(A)$ with $A \subseteq E$ and
   $A E \subseteq E$ (respectively, $E A \subseteq E$) and let $\var: E
   \ra I(A)$ be a completely bounded left (resp., right) $A$-module
   map. If $\var (a) = a$ for all $a$ in $A,$ then $\var (e) = e$ for all
   $e$ in $E.$
\end{cor}

\begin{proof}
There exists $y$ in $I(A)$ with $\var (e)  = e \cdot y$
for all $e$ in $E.$  Since $\var (1) = 1, \  y = 1$ and hence, $\var (e) = e$
for all $e$ in $E.$ 
\end{proof}

\section{Injective Bimodule Extensions of C*-algebras $A$ and the Injective
Envelope $I(A)$}

Let $A$ and $B$ be unital $C^*$-algebras. Recall the definition of {\it
operator $A-B$-bimodules}. These are operator
spaces $E$ which are $A-B$-bimodules and such that the trilinear module
pairing $A \times E \times B \ra E, \ (a, e, b) \ra aeb$ is completely
contractive in the sense of E.~Christensen and A.~Sinclair \cite{C-S}.
This is equivalent to requiring that for matrices $(a_{i,j}), (e_{i,j})$
and $(b_{i,j})$ of the appropriate sizes, the induced matricial module
product is contractive, i.e.,
$$
  \left\| \left( \sum_{k,m} a_{i,k}e_{k,m}b_{m,j} \right) \right\| \le
  \|(a_{i,j})\| \|(e_{i,j})\|\|(b_{i,j})\| \, .
$$
These are the objects of the category ${}_A\cO_B,$ \cite{BMP}, \cite{P2} and
the morphisms between two operator $A-B$-bimodules in this category are the
completely bounded $A-B$-bimodules maps.
When we want to restrict the morphisms to be completely contractive
$A-B$-bimodule maps we will denote
the category by ${}_A\cO^1_B$.

We assume that all module actions are unital, i.e.
$1\cdot e \cdot1 = e.$  We set ${}_A{\cO} = {}_A{\cO}_\C$ and call these
{\it left operator $A$-modules}, and $\cO_A = {}_\C \cO_A$ and call
these {\it right operator $A$-modules}.

\begin{dfn}
  An operator $A-B$-bimodule $I$ is {\it $A-B$-injective}, if
  whenever $E \subseteq F$ are operator $A-B$-bimodules, then every
  completely bounded $A-B$-bimodule map from $E $ into $I$ has a completely
  bounded $A-B$-bimodule extension to $F.$  Note that we do {\it not}
  require that the $cb$-norm of the extension is the same as the $cb$-norm
  of the original map. When this is the case we will call $I$ a {\it tight
  $A-B$-injective $A-B$-bimodule}.
\end{dfn}

Some comments on terminology are helpful.
Our definition of $A-B$-injective is the usual definition of injectivity in
the category ${}_A\cO_B$, while what we are calling tight $A-B$-injective is
the corresponding definition of injectivity in
the category ${}_A\cO^1_B$.

If $A$ and $B$ are both
$C^*$-subalgebras of $B(H),$ then by the bimodule version of G.~Wittstock's
extension theorem \cite[Thm.~4.1]{Wi} (see also \cite{Suen}) $B(H)$ is a
tight
$A-B$-injective. Thus, if $M \subseteq B(H)$ is the range of a completely
bounded projection $\var : B(H) \ra M,$ which is also an $A-B$-bimodule map,
then $M$ is $A-B$-injective, but it is not evidently tight $A-B$-injective.
A $C^*$-subalgebra $I \subseteq B(H)$
is generally called ``injective" if it is the range of a completely positive
projection.  This term is so widespread we continue to use it here. Note
that such a map is also automatically an $I$-bimodule map.  Thus,
such an $I$ is a
tight $A-B$-injective for any $C^*$-subalgebras $A$ and $B$ of $I.$

In particular, M.~Hamana's injective envelope $I(A)$ is a tight
$A-A$-injective $A-A$-bimodule,
a tight $A-\C$-injective left $A$-module and a tight $\C-A$-injective right
$A$-module.

On the other hand there are many $\C-\C$-injectives which are not tight,
i.e.~not injective in the usual sense. For example, for any subspace $E$
of $B(H)$ of finite codimension, it is easy to show that there
exists a completely bounded projection from $B(H)$ onto $E$ and hence $E$ is
$\C-\C$-injective.

T.~Huruya \cite{Hu} has given an example of a unital $C^*$-subalgebra of an
injective $C^*$-algebra of finite codimension that is not injective.
By the above argument, this algebra is the range of a completely bounded
projection and hence is $\C-\C$-injective.
Thus there exist $C^*$-algebras that are $\C-\C$-injective, but are not
injective in the usual sense.

In our terminology, M.~Hamana's rigidity result implies that
if $A \subseteq E \subseteq I(A)$ and $E$ is a tight $\C-\C$-injective,
then $E = I(A).$ We prove this fact in the remarks following the theorem.

\begin{thm}  \label{thm-2}
   Let $A$ be a unital $C^*$-algebra and let $A \subseteq E \subseteq I(A).$
Then the following are equivalent:
     \begin{itemize}
      \item [a)] $E$ is $A-\C$-injective,
      \item [b)] $E$ is $\C-A$-injective,
      \item [c)] $E$ is $A-A$-injective,
      \item [d)] $E = I(A).$
     \end{itemize}
\end{thm}

\begin{proof}
By the Hahn-Banach extension theorem for completely bounded $A-B$-bimodule
maps \cite[Thm.~4.1]{Wi} (see also \cite{Suen}, \cite{Pa}), it follows that
$I(A)$ is $A-\C$-injective, $\C-A$-injective and $A-A$-injective. Thus, d)
implies a), b) and c).
It will suffice to prove that a) implies d), the other implications are
similar.

If $E$ is $A-\C$-injective then the identity map from $E$ to $E$ extends to
a
completely bounded left $A$-module projection from $I(A)$ to $E.$  Letting
$I(A)$ play the role of $E$ in the rigidity theorem yields the result.
\end{proof}

The module actions are necessary in the above theorem. Since there always
exists a completely bounded projection from any operator space onto
a subspace of finite codimension, if $A \subseteq E
\subseteq I(A)$ with $E$ a subspace of finite codimension then $E$ is
$\C-\C$-injective but $E \neq I(A).$

On the other hand, if we required $E$
to be tight, then there would exist a completely contractive projection
$\var$ onto $E.$  Since 1 belongs to $E$ we would have $\var (1) = 1$ and,
consequently,
this projection would be completely positive.  Hence $E$ would be an
operator system.  Thus, $E = I(A)$ by M.~Hamana's rigidity theorem and we
would be adding nothing new.

We now are in a position to clarify the relationship between these new
notions of injectivity and injectivity in the usual sense for C*-algebras.

\begin{thm}  \label{thm-inj}
     Let $A$ be a unital C*-algebra. Then the following are equivalent:
       \begin{itemize}
         \item [a)] $A$ is an injective C*-algebra(in the usual sense),
         \item [b)] $A$ is a tight $A-\C$-injective module,
         \item [c)] $A$ is a $A-\C$-injective module,
         \item [d)] $A$ is a tight $\C-A$-injective module,
         \item [e)] $A$ is a $\C-A$-injective module,
         \item [f)] $A$ is a tight $A-A$-injective module,
         \item [g)] $A$ is a $A-A$-injective module.
          \end{itemize}
\end{thm}

\begin{proof}
We prove the equivalence of a),b) and c), the remaining arguments are
similar. We have that a) implies b) by Wittstock's Hahn-Banach extension
theorem
for module maps. Clearly, b) implies c). We now prove that c) implies a).
Since $A$ is a $A-\C$-injective module, the identity map on $A$ extends to
a completely bounded left $A$-module map from $I(A)$ into $A$.
But by the Rigidity Theorem, this extended map must be the identity map on
$I(A)$ and hence $I(A) = A$. Thus, $A$ is injective.
\end{proof}
\begin{dfn}
   Let $M$ be an operator $A-B$-bimodule. Call $I$ a {\it minimal
   $A-B$-injective extension} of $M$, if $M \subseteq I$ and $M \subseteq E
   \subseteq I$ with $E \ A-B$-injective implies $E = I.$

   By Theorem \ref{thm-2}, $I(A)$ is a minimal $A-\C$-injective extension of
   $A$ and also a minimal $\C-A$ and $A-A$ injective extension of $A.$

   We call a map $\var$ a {\it completely bounded isomorphism} if
   both $\var$ and $\var^{-1}$ are completely bounded.
\end{dfn}

\begin{thm} \label{thm-3}
   Let $A$ be a unital $C^*$-algebra and let
   $I$ be a minimal $A-\C$-injective extension of $A,$ then there exists
   a completely bounded left $A$-module isomorphism $\var : I(A) \ra I$ with
   $\var (a) = a$ for all $a$ in $A.$  If we require that $I$ is also a
tight
   $A-\C$-injective then $\var$ may be taken to be a complete isometry.
   Analogous statements hold for right modules and bimodules.
\end{thm}

\begin{proof}
Since $I$ and $I(A)$ are $A-\C$-injective, there exist completely
bounded left $A$-module maps $\var : I(A) \ra I$ and $\psi : I \ra I(A)$
which fix $A.$  By the rigidity of $I(A), \ \psi \circ \var$ is the identity
on $I(A)$ and hence $\var \circ \psi$ is the identity restricted to $E = $
range $(\var).$ This makes $E$ an $A-\C$-injective module and
hence $E = I$ and $\psi = \var^{-1}.$
\end{proof}

If we required $I$ to be tight and only minimal among all tight injectives,
then as in the remark following Theorem \ref{thm-2}, our result would reduce
to M.~Hamana's theory.

We now turn to some applications to projections. In \cite{Bu-Pa} it was
shown
that if $M \subseteq B(H)$ is a von Neumann algebra and there exists a
bounded
$M$-bimodule projection, $\var : B(H) \ra M,$ then $M$ is injective. Such a
map $\var$ is easily shown to be automatically completely bounded.

In \cite{P2} the same result was shown to hold for $C^*$-algebras. The above
results on injective envelopes allow us to extend these results a bit.
Perhaps, more importantly, the new proof is much simpler than the proof in
\cite{P2}.

\begin{thm}  \label{thm-4}
   Let $A \subseteq B(H)$ be a unital $C^*$-algebra. If there exists
   a completely bounded left (or right) $A$-module projection of $B(H)$ onto
   $A,$ then $A$ is injective.
\end{thm}

\begin{proof}
Since $B(H)$ is $A-A$-injective the identity map on $A$ extends to
a completely bounded $A$-bimodule map from $I(A)$ to $B(H).$  Composing with
the projection onto $A$ gives a completely bounded left $A$-module map from
$I(A)$ to $A$ which is the identity on $A.$  By rigidity (Corollary
\ref{cor-rigid}) $I(A) = A$ and hence $A$ is injective.
\end{proof}

The following example gives an indication of the obstacles that arise in
attempting to generalize Theorem \ref{thm-1},
Corollary \ref{cor-rigid} and Theorem \ref{thm-4} to the
case of non-involutive operator algebras in $B(H).$
Consider the operator algebra $A
\subset M_2({\mathbb C})$ defined by
\[
  A = \left\{ X \in M_2({\mathbb C}) \, : \, X = S^{-1} \cdot {\rm
diag}(a,b)
  \cdot S \, , \, a,b \in {\mathbb C} \right\} \, , \,\,\,
  S = \left( \begin{array}{cc} 1 & 1 \\ 0 & 1 \end{array} \right) \, .
\]
If $\omega: M_2({\mathbb C}) \to {\mathbb C} \oplus {\mathbb C} \subset
M_2({\mathbb C})$ is the canonical conditional expectation on $M_2({\mathbb
C})$
preserving the diagonal and mapping off-diagonal elements to zero, then the
map $\phi: M_2({\mathbb C}) \to A$ defined by the rule $\phi(X) = S^{-1}
\cdot
\omega(SXS^{-1}) \cdot S$ is a completely bounded $A$-bimodule projection.
However, the smallest $C^*$-subalgebra of $M_2({\mathbb C})$ generated by
$A$
is $M_2({\mathbb C})$ itself, and since $M_2({\mathbb C})$ is injective we
obtain $I(A + A^*) = M_2({\mathbb C})$.
Consequently, Theorem \ref{thm-1} and Corollary \ref{cor-rigid} cannot be
extended to this situation, and Theorem \ref{thm-4} is not true for the
described operator algebra $A$.

\smallskip
In \cite{C-S2} and \cite{Pi2} it was proven that if $M \subseteq B(H)$ is
a von Neumann algebra and if there exists a completely bounded projection of
$B(H)$ onto $M$ then $M$ is injective, cf.~\cite{Pi3}. The direct analogue of
this result is false for $C^*$-algebras. T.~Huruya \cite{Hu}
gave an example of a
non-injective $C^*$-subalgebra of codimension 1 of an injective
$C^*$-algebra.
It is easily shown that any time Huruya's algebra is represented as a
$C^*$-subalgebra of $B(H),$ then there will exist a completely bounded
projection of $B(H)$ onto this non-injective algebra.
Thus, to generalize the results of \cite{C-S2} and \cite{Pi2}, we will need
an additional condition.
\begin{dfn}
   An operator $A-B$-bimodule $R$ is {\it relatively $A-B$-injective}
   if whenever $E \subseteq F$ are operator $A-B$-bimodules such that there
   exists a completely bounded projection of $F$ onto $E,$ then every
   completely bounded $A-B$-bimodule map of $E$ into $R$ has a completely
   bounded $A-B$-bimodule extension to $F.$ It is important to note that we
   do not require that the projection of $F$ onto $E$ is an $A-B$-bimodule
   map.
\end{dfn}

The concept of relative injectivity was introduced
 in \cite{P2} with slightly
different notation, relative $A-B$-injective was denoted $(A-B,
\C-\C)$-injective. A $C^*$-algebra $A$ is $A-A$-injective if and only if $A$
is injective in the usual sense. In contrast, \cite{P2} showed that
every von Neumann algebra $M$ is relatively $M-M$-injective,
$M-\C$-injective and $\C-M$-injective.

\begin{thm}  \label{thm-5}
   Let $A\subseteq B(H)$ be a unital $C^*$-subalgebra. If there exists a
   completely bounded projection of $B(H)$ onto $A$ and $A$ is relatively
   $A-\C$-injective (or $\C-A$-injective), then $A$ is injective.
\end{thm}

\begin{proof}
Since $A$ is $A-\C$-injective, the identity map from $A$ to $A$
has a completely bounded left $A$-module extension to $B(H).$ This map is
clearly a projection. Hence $A$ is injective by Theorem \ref{thm-4}.
\end{proof}

Because von Neumann algebras are relatively injective \cite{P2}, Theorem 3.5
implies the result of \cite{C-S2} and \cite{Pi2}.

\begin{cor}
   Let $M \subseteq B(H)$ be a von Neumann algebra. If there exists a
   completely bounded projection of $B(H)$ onto $M,$ then $M$ is injective.
\end{cor}

Relative injectivity was shown in \cite{P2} to be equivalent to the
vanishing of
certain completely bounded "Ext" groups, which in turn implied the vanishing
of completely
bounded Hochschild cohomology. Thus, relative injectivity captures both the
vanishing of cohomology
and these projection results. It is still unknown which
$C^*$-algebras $A$ are relatively $A-\C$-injective.
By Theorem \ref{thm-5}, T.~Huruya's $C^*$-algebra $A,$ cannot be relatively
$A-\C$-injective.

\section{Local Multiplier Algebras, Injective Envelopes and Regular
             Completions}

We close this paper with some applications to multiplier algebras.
Our main point is that by invoking Theorem \ref{thm-1}, we will
see immediately that multipliers are "naturally" represented as
multiplication by elements in $I(A).$ This concrete representation of
multipliers can be used to simplify some arguments. Thus, Theorem
\ref{thm-1} provides an
alternative starting point for developing the theory of multipliers.

In particular, we show
that $I(A)$ contains the local multiplier algebra of $A, \ M_{loc} (A),$
intrinsically as a $C^*$-subalgebra. Recall that a closed 2-sided ideal $J$
of $A$ is called {\it essential} if $J \cap K\neq \{ 0 \}$ for every
non-trivial 2-sided ideal $K.$ All ideals in this section are norm-closed.

The {\it left multiplier algebra $LM(J)$ of $J$} is just the set of right
$A$-module maps $\psi : J \ra J.$ Such a map is automatically (completely)
bounded and $\V{\psi} = \V{\psi}_{cb}.$ The {\it right multiplier algebra}
$RM(J)$ is defined similarly. The {\it multiplier algebra} $M(J)$ consists
of pairs of linear maps $\var, \psi :J \ra J$ satisfying $\var (j_1) j_2 =j_1
\psi (j_2).$ This identity implies that $\psi \in LM(J), \ \var \in RM(J)$
and $\V{\var} = \V{\psi} = \V{\var}_{cb} = \V{\psi}_{cb}.$

The {\it local multiplier algebra} $M_{loc} (A)$ is defined by taking a
direct limit of $M(J)$ over all essential ideals $J$ of $A$ ordered by
reverse inclusion. See \cite{AM} for details.

\begin{lem}  \label{lem-1}
   Let $A$ be a unital $C^*$-algebra and let $J$ be a 2-sided essential
   ideal of $A$ and let $\var : J \ra I(A)$
   be a completely bounded left (resp., right) $A$-module map. Then there
   exists a unique element $x$ in $I(A)$ such that $\var$ is right (resp.,
   left) multiplication by $x.$ Moreover, $\V{x} = \V{\var} =
   \V{\var}_{cb}.$
\end{lem}

\begin{proof}
By Theorem \ref{thm-1} such an $x$ exists, it remains to show that $x$
is unique. To this end consider, $F = \{y \in I(A): Jy = 0\}$ which is
clearly a right $A$-submodule of $I(A).$ It will suffice to show that
$F = \{ 0 \}.$ Let $\{e_{\al}\}$ be a contractive approximate identity for
$J.$ For $a \in A, \ y \in F,$ we have,
  $$
    \V{a-y} \geq \sup_{\al} \V{e_{\al}(a-y)} = \sup_{\al} \V{e_{\al} a} =
\V{a}
  $$
with the last equality using the fact that $J$ is essential. The same
calculation for matrices shows that the quotient map $q: I(A) \ra I(A)/F$ is
a complete isometry on $A$ and a right $A$-module map. Now
since $I(A)$ is injective, there exists a
completely contractive right $A$-module map $\var : I(A)/F \ra I(A)$.  Hence
by rigidity $\var \circ q (b) = b$ for all $b$ in $I(A),$ and it follows
that
$F = \{ 0 \}.$
\end{proof}

The fact that $F$ must be $\{ 0 \}$ is related to the fact that $I(A)$
is in a certain sense an ``essential extension" of $A.$

\begin{thm}  \label{thm-6}
   Let $A$ be a unital $C^*$-algebra, let $J$ be a 2-sided essential ideal in 
   $A$ and let $(\var, \psi)$ be in $M(J).$ Then there exists a unique element 
   $x$ in $I(A)$ such that $\var (j) = jx, \psi (j) = xj$ for all $j$ in $J.$
\end{thm}

\begin{proof}
By Lemma \ref{lem-1}, there exist unique elements $x_1, x_2$ in $J$ such
that $\var (j_1) = j_1 x_1, \ \psi (j_2) = x_2 j_2$ for all $j_1, J_2$ in
$J.$  But $j, x_1 j_2 = \var (j_1) j_2 = j_1 x (j_2) = j_1 x_2 j_2$ and so
$j_1 (x_1 - x_2)j_2 = 0$ for all $j_1.$  Applying Lemma \ref{lem-1} we
conclude that $(x_1 - x_2)j_2 = 0$ for all $j_2$ and so $x_1 = x_2.$
\end{proof}

\begin{cor}  \label{cor-2}
   The inclusion of $A$ into $I(A)$ extends in a unique way to a
   $*$-monomorphism of $M_{loc} (A)$ into $I(A).$ The image of $M_{loc}(A)$
   under this map is
     $$
       \{x \in I(A): xJ \subseteq J \text{ and } Jx \subseteq J
       \text{ for some essential ideal $J$}\}^{--}.
     $$
\end{cor}

\begin{proof}
For each $(\var, \psi)$ in $M(J)$ there exists a unique $x$ in $I(A)$
implementing $(\var, \psi).$ By this uniqueness the map $(\var, \psi) \ra x$
must be a $*$-monomorphism on $M(J).$  Furthermore, let $J_i$ be essential
ideals, and let $(\var_i, \psi_i)$ in $M(J_i)$ be implement by $x_i.$  If
$\var_1 = \var_2$ and $\psi_1 = \psi_2$ on $J_1 \cap J_2$ then, since $J_1
\cap J_2$ is essential, we must have $x_1 = x_2.$

This shows that the inclusions of $M(J_i)$ into $I(A)$ are coherent and
allows us to extend these $*$-monomorphism to the direct limit, $M_{loc}
(A).$

Now assume that $\pi : M_{loc} (A) \ra I(A)$ is any $*$-monomorphism with
$\pi (a) = a$ for all $a$ in $A.$  Then for $(\var, \psi)$ in $M(J)$
  $$
    \pi ((\var, \psi)) j = \pi ((\var, \psi) j) = \pi (\psi (j)) = \psi (j),
  $$
and $j \pi ((\var, \psi)) = \var(j)$ from which it follows that $\pi ((\var,
\psi))$ is the unique element implementing $(\var, \psi).$

Finally, since $\{x \in I(A): xJ \subseteq J \text{ and } Jx \subseteq J\}$
is exactly the image of $M(J)$ we have the last claim.
\end{proof}

\begin{rmk}
   The above results allow one to define a {\it local left} (resp.,
   {\it right}) {\it multiplier algebra of $A$} easily,
   which we have not seen in the literature. Indeed, if $LM_{loc} (A) =
   \{x \in I(A): xJ \subseteq J \text{ for some essential ideal
   $J$}\}^{--},$
   then this set is easily seen to be completely isometrically
   isomorphic to the direct limit of $LM(J).$ It is interesting to note that
   if $J_1$ and $J_2$ are essential ideals and $\var \in LM(J_1)$ then
   $\var (J_1 \cap J_2) \subseteq J_1 \cap J_2$ and by Lemma \ref{lem-1},
   $\V{\var} = \V{\var\mid_{J_1 \cap J_2}}.$ We
   define the local right multiplier algebra $RM_{loc}(A)$ of $A$, analogously.

   To define the {\it local quasi-multiplier space $QM_{loc}(A)$ of $A$} we
   have to recall that the injective envelope $I(A)$ of $A$ is a monotone
   complete $C^*$-algebra and, hence, an $AW^*$-algebra. On the other
   hand norm-closed two-sided ideals $J$ of $C^*$-algebras are
   automatically hereditary, and so we can apply \cite[Cor.~1.4]{Fr99}:
   for every norm-closed two-sided ideal $J \subseteq A$ and every
   quasi-multiplier $x \in QM(J)$ there exists an element
   $\overline{x} \in I(A)$ such that $j_1xj_2=j_1\overline{x}j_2$ for
   any $j_1,j_2 \in J$ and $\| \overline{x} \|$ equals the norm of $x$
   estimated in the bidual von Neumann algebra $A^{**}$. For essential
   ideals $J$ the element $\overline{x}$ has to be unique, in fact it
   can be found as a quasi-strict limit of nets of $J$. Since inclusion
   relations of essential ideals and their corresponding quasi-multiplier
   spaces are respected inside $I(A)$ we can define $QM_{loc} (A) = \{x
   \in I(A): JxJ \subseteq J \text{ for some essential ideal $J$}\}^{--},$
   to be the local quasi-multiplier space of $A$.
\end{rmk}

Note that in any situation where $M_{loc} (A) \not\equiv LM_{loc} (A)$
then necessarily $M_{loc}(A) \not\equiv I(A).$ (In general, the conditions
$M_{loc}(A) \not\equiv LM_{loc}(A)$, $M_{loc}(A) \not\equiv QM_{loc}(A)$
and $LM_{loc}(A) \not\equiv QM_{loc}(A)$ are equivalent by general
multiplier theory.) If $A$ is any simple, unital, non-injective
$C^*$-algebra
like a non-injective type ${\rm II}_1$ or type III von Neumann factor then
$A = M_{loc} (A) = LM_{loc}(A) \neq I(A).$ However, in some cases we obtain
the coincidence of the $C^*$-algebras $M_{loc}(A) = I(A).$

\begin{prop}
   Let $A$ be a unital $C^*$-algebra, $K \subseteq A \subseteq B(H)$ where
   $K$ denotes the ideal of compact operators, then $M_{loc}(A) = I(A) =
   B(H),$ $*$-isomorphically.
\end{prop}

\begin{proof}
Since $K$ is necessarily an essential ideal of $A$ and $M(K) = B(H)$
we have $B(H) \subseteq M_{loc}(A).$ By Corollary \ref{cor-2} we have a
$*$-monomorphism $\pi$ of $M_{loc} (A)$ into $I(A).$ Hence, $A \subseteq
B(H) \subseteq M_{loc}(A) \subseteq I(A),$ as $C^*$-algebras.
Since $B(H)$ is $A-A$-injective, by Theorem 3.1, we have $B(H) = I(A)$ and
the result follows.
\end{proof}

The fact that $I(A) = B(H)$ is due to M.~Hamana \cite{Ham} with a different
proof.

\begin{thm}  \label{thm-7}
   Let $A$ be a commutative unital $C^*$-algebra, then $M_{loc} (A) = I(A),$
   $*$-isomorphically.
\end{thm}

\begin{proof}
By \cite[Thm.~1]{AM} $M_{loc} (A)$ is a commutative $AW^*$-algebra. However,
commutative $AW^*$-algebras are injective by \cite[Th.~25.5.1]{Sem} since
bounded linear maps between C*-algebras are positive whenever their norm
equals their evaluation at the identity of the C*-algebra. Consequently, the
$*$-monomorphism of $M_{loc} (A)$ into $I(A)$ must be onto.
\end{proof}

In the theory of local multiplier $C^*$-algebras the problem whether
$M_{loc}(A)$ coincides with $M_{loc}(M_{loc}(A))$ for any $C^*$-algebra $A$
is one of the main open questions, cf.~\cite{AM,So2}.
Set $M_{loc}^{k+1}(A) = M_{loc}(M_{loc}^{k}(A))$, which is called the
$(k+1)$-order local multiplier algebra of $A$.
We show that any higher order local multiplier $C^*$-algebra of a given
$C^*$-algebra $A$ is contained in its injective envelope $I(A)$ and, what
is more, that the injective envelopes $I(A)$ and $I(M_{loc}^{k}(A))$ coincide
for any $C^*$-algebra $A$. The latter is of
special interest since general $C^*$-subalgebras $A$ of injective
$C^*$-algebras $B$ might not admit an embedding of their injective envelopes
$I(A)$ as a $C^*$-subalgebra of $B$ that extends the given embedding of $A$
into $B$, see \cite[Rem.~3.9]{Ham2} for an example.

\begin{thm}
   Let $A$ be a unital $C^*$-algebra and $M_{loc}(A)$ be its local multiplier
   $C^*$-algebra. Then the injective envelope $I(A)$ of $A$ is the
   injective envelope $I(M_{loc}(A))$ of $M_{loc}(A)$ and consequently,
   $M_{loc}^{k}(A)$ is contained in $I(A)$ for all $k$.
\end{thm}

\begin{proof}
Since $M_{loc}(A)$ is $*$-isomorphically  embedded into $I(A)$ extending the
canonical $*$-monomorphism of $A$ into $I(A)$ by Theorem \ref{thm-1}, the
$C^*$-algebra $I(A)$ serves as an injective extension of the $C^*$-algebra
$M_{loc}(A)$, cf.~\cite{Ham}.
However, the identity map on $M_{loc}(A)$ admits a unique extension to a
completely positive map of $I(A)$ into itself with the same completely
bounded norm one since $A \subseteq M_{loc}(A) \subseteq I(A)$ by
construction and $I(A)$ is the injective envelope of $A$. So $I(A)$ has
to be the injective envelope of $M_{loc}(A)$, too.
\end{proof}

\begin{prob}
   Characterize the $C^*$-algebras $A$ for which the local
   multiplier $C^*$-algebra $M_{loc} (A)$ of $A$ coincides with the
   injective envelope $I(A)$ of $A$ or at least with the regular monotone 
   completion $\overline{A}$ of $A$ in $I(A)$.
\end{prob}

This question is surely difficult to answer: if $A$ is an $AW^*$-algebra
then
the local multiplier algebra $M_{loc}(A)$ of $A$ coincides with $A$ by
\cite{Ped}. However, $A$ coincides with its regular monotone completion
$\overline{A}$ if and only if $A$ is monotone complete. So we arrive at a
long standing open problem of $C^*$-theory dating back to the work of
I.~Kaplansky in 1951 (\cite{Kapl}): Are all $AW^*$-algebras monotone
complete, or do there exist counterexamples?

\begin{rmk}
  If $A$ is a non-unital $C^*$-algebra and $B$ denotes its unitization, then
  $A$ is a 2-sided essential ideal in $B$.
  Hence, by Theorem \ref{thm-6}, $M(A) \subseteq I(B).$
  However, in \cite{BP}, it is observed that $I(A) = I(B)$, and so the hypothesis 
  that $A$ is unital can be removed from Theorem \ref{thm-6}.
  Similarly, every 2-sided essential ideal in $A$ is an essential ideal in $B$, so 
  that Corollary \ref{cor-2} applies for non-unital $A$ as well.
  Similar arguments show that the unital hypothesis can be dropped in Proposition 
  4.4, Theorem 4.5 and Theorem 4.6.
\end{rmk}


\end{document}